\documentclass[10pt]{article}
\usepackage{amsmath}
\usepackage{amsfonts}
\usepackage{amssymb}

\begin{document}

\numberwithin{equation}{section}

\def\ye{y_{\varepsilon}}

\begin{center}
{\bf BOUNDARY VALUE PROBLEM FOR A MULTIDIMENSIONAL SYSTEM OF EQUATION WITH 
RIEMANN-LIOUVILLE
DERIVATIVES}
\end{center}
\begin{center}
{\bf M.O. Mamchuev }
\end{center}

In the paper а boundary-value problem for a multidimensional system of partial differential 
equations with fractional derivatives in Riemann-Liouville 
sense with constant coefficients 
is studied in a rectangular domain. The existence and uniqueness theorem for the solution 
of the boundary value problem is proved.
The solution is constructed in explicit form in terms of the Wright function of the 
matrix argument.

 \medskip


 \smallskip


\section{Introduction}\label{sec:1}

\setcounter{section}{1}
\setcounter{equation}{0}

Consider the system of equations
\begin{equation} \label{bsn}
\sum\limits_{i=1}^{m}A_i D_{0x_i}^{\alpha_i} u(x)=Bu(x)+f(x),
\quad 0<\alpha_i<1,  
\end{equation}
in the domain $\Omega=\{x=(x_1,...,x_m):\ 0<x_i<a_i\leq \infty, \, i=\overline{1,m}\},$ 
where  
$f(x)=||f_1(x),...,f_n(x)||$ and 
$u(x)=||u_1(x),...,u_n(x)||$ are 
the given and required $n$-dimensional vectors, respectively, 
$A_i,$ $B$  are given constant square matrices of order $n,$
$D_{ay}^{\nu}$ 
is the operator of fractional integro-differentiation in Riemann-Liouville sense 
of order $\nu.$ 

Operator $D_{ay}^{\nu}$ is determined 
for $\nu< 0$ by the following formula \cite[p. 9]{Nakhushev-2003}: 
$$D_{ay}^{\nu}g(y)=\frac{\mathop{\rm sgn}(y-a)}{\Gamma(-\nu)}
\int\limits_a^y \frac{g(s)ds}{|y-s|^{\nu+1}}, $$
and for $\nu \geq 0$ 
can be determined with the help of following recursive relation
$$D_{ay}^{\nu}g(y)=\mathop{\rm sign}(y-a)\frac{d}{dy}D_{ay}^{\nu-1}g(y),$$
$\Gamma(z)$ is a Eullers gamma-function.

Let all the eigenvalues of matrices $A_i$ $i=\overline{1,m}$ are positive. 
Without loss of generality, we assume that $A_1 = I$ is the identity matrix of order $ n. $

We formulate the boundary value problem for the system (\ref{bsn}).

{\bf Problem 1.} 
{\it  Find the solution $u(x)$ of the system (\ref{bsn}) satisfying the following boundary conditions
\begin{equation} \label{bu2}
\lim\limits_{x_i \rightarrow 0} D_{0x_i}^{\alpha_i-1} u=\varphi_i (x_{(i)}), \quad
x_{(i)}\in\Omega^i,  \quad i=\overline{1,m}, 
\end{equation}
where $\!x_{(i)}\!=\!(x_1,...,x_{i-1},x_{i+1},...,x_m), 
\Omega^i\!=\!\omega_{a_1}\times...\times\omega_{a_{i-1}}\times\omega_{a_{i+1}}\times...
\times\omega_{a_m},$
$\omega_{a_j}=\{x_j:0<x_j<a_j\},$
$\varphi_i (x_{(i)})$ are given $n$-dimensional vectors. }

Let us review the works  associated with the investigation of the system 
(\ref{bsn}) including in the scalar case $n=1.$
In paper \cite{Clement-2000} 
for the equation
\begin{equation} \label{cgl}
D_{0x}^{\alpha} (u-h_1(y))+D_{0y}^{\beta} (u-h_2(x))=f,   \quad 0<\alpha, \beta<1, \quad x,y\geq 0,
\end{equation}
the solvability of the boundary value problem is studied in a class of 
Helder's continuous functions 
with the initial conditions
$u(0,y)=h_1(y),$   $u(x,0)=h_2(x)$ 
and the right hand side
$f(x,y).$
The fundamental solution of the equation  (\ref{cgl}) was rewriten in the form
\begin{equation} \label{fscgl}
\psi_{\alpha,\beta} (x,y)= \int\limits_0^{\infty}\tau^{-\frac{1}{\alpha}-\frac{1}{\beta}}
\varphi_{\alpha}\left(x\tau^{-\frac{1}{\alpha}}\right)\varphi_{\beta}\left(y\tau^{-\frac{1}{\beta}}\right)d\tau,
\end{equation}
where 
\[\varphi_{\mu}(t)=\frac{1}{\pi}\sum\limits_{k=1}^{\infty}\frac{(-1)^{k+1}}{k!}
\sin (\pi\mu k)\Gamma(\mu k+1)t^{-\mu k-1}.\]
Using the equality 
$\Gamma(1-z)\Gamma(z)\sin(\pi z)=\pi,$
the function $\varphi_{\mu}(t)$ can be represented in the form
\[\varphi_{\mu}(t)=\frac{1}{t}\phi(-\mu,0;-t^{-\mu}),\]
where $\phi(\rho,\mu;z)=\sum\limits_{k=0}^{\infty}\frac{z^k}{k!\Gamma(\rho k+\mu)}$ is the Wright's function
\cite{Wright-1933}, 
\cite{Wright-1934}. 
From the last relation and (\ref{fscgl}) we get the following representation
\[\psi_{\alpha,\beta} (x,y)= \frac{1}{xy}\int\limits_0^{\infty}
\phi(-\alpha,0; -\tau x^{-\alpha})\phi(-\beta,0;-\tau y^{-\beta})d\tau.
\] 
Holder's smuthness of the following equation's solution
$$ D_{0x}^{\alpha}(u-u_0)+c(x,y)u_y(t,x)=f(x,y), \quad x,y>0,$$
satisfaing the boundary conditions  $u(0,y)=u_0(y)$ and $u(x,0)=u_1(x),$
was studied in paper \cite{Clement-1999}.

In the papers
\cite{Pskhu-2003}, 
\cite{PsMon}, 
the theorems of uniquiness and existence of regular solution are prooved 
for a boundary value problem in a rectangular domain for the equation
\begin{equation} \label{ps1}
D_{0x}^{\alpha} u(x,y)+\lambda D_{0y}^{\beta} u(x,y)+\mu u(x,y)=f(x,y),   
\quad 0<\alpha, \beta<1, \, \lambda>0, \, x,y> 0.
\end{equation}
When $\lambda=1$ the fundamental solution has the form
\[w(x,y)= \frac{1}{xy}\int\limits_0^{\infty} e^{-\mu \tau}
\phi(-\alpha,0; -\tau x^{-\alpha})\phi(-\beta,0;-\tau y^{-\beta})d\tau,
\] 
and when $\mu=0$  it has the form
\[w(x,y)=\frac{x^{\alpha-1}}{y}e_{\alpha,\beta}^{\alpha,0}\left(-\lambda\frac{x^{\alpha}}{y^{\beta}}\right),\]
where
\[e_{\alpha,\beta}^{\mu,\nu}(z)=\sum\limits_{k=0}^{\infty}\frac{z^k}{\Gamma(\mu+\alpha k)\Gamma(\nu-\beta k)}\]
is the Wrihgt type function \cite{PsMon}.
In the case $\alpha=1,$ $\mu=0,$ a boundary value problem with negative coefficient $\lambda<0$
was studied for the equation (\ref{ps1}).

The equation (\ref{ps1}) 
with  variable coefficients $\lambda\equiv \lambda(x)$ and $\mu\equiv\mu(x),$ 
when $\alpha=1$ and $\lambda(x)$ can have a zero of order $m\geq 0$ at the point $x=0$ 
was investigated in the papers
\cite{Mamchuev-2009-1}, 
\cite{Mamchuev-2009-2}, 
\cite{MamMon}. 
It's fundamental solution  
\[w(x,y;t,s)=\frac{\exp(\Lambda(x,t))}{y-s}\phi\left(-\beta,0;-M(x,t)(y-s)^{-\beta}\right),\]
where $\Lambda(x,t)=\int\limits_t^x\lambda(\xi)d\xi,$ $M(x,t)=\int\limits_t^x\mu(\xi)d\xi,$
was constracted, and  
the existence and uniqueness theorems for the solutions of the boundary-value problem in a rectangular domain and the Cauchy problem were proved.

In \cite{Pskhu-2011} 
the unique solvability of the analogue of problem 1 for the equation (\ref{ps1}) 
with the Dzhrbashyan-Nersesyan fractional differentiation operators  was investigated
in the case
$n=1,$ $A_i=\lambda_i>0$ $(i=\overline{1,m}),$ $B=\lambda_0$,
and the fundamental solution
$$w(x)=\int\limits_0^{\infty}e^{-\lambda_0\tau}
\prod\limits_{i=1}^{m}
\frac{1}{x_i}\phi\left(-\alpha_i,0;-\lambda_i\tau x_i^{-\alpha_i}\right)d\tau$$
of this equation was constracted. 

For two-dimensional system
\begin{equation} \label{sn}
D_{0x}^{\alpha} u(x,y)+AD_{0y}^{\beta} u(x,y)=Bu(x,y)+f(x,y),
\end{equation}
the problem 1 was solved in explicite form in 
\cite{Mamchuev-2008}, 
when $A$ was an identity matrix and in
\cite{Mamchuev-2017},
when $A$ was positive defined matrix.
In the paper \cite{Mamchuev-2017} the fundamental solution of the system (\ref{sn})
was constracted in the terms of the introdused Wright's function of
the matrix argument.

In present paper we will use similar approuch for the solving
of problem 1 in multidimensional case.

\section{Wright's function}\label{sec:2}

\setcounter{section}{2}
\setcounter{equation}{0}

Wright function  
\cite{Wright-1933}, 
\cite{Wright-1934} 
is called an entiere function, which is depended from two parameters $\rho$ and $\mu,$ 
and representad by the series
\begin{equation*}
 \phi(\rho,\mu;z)=\sum\limits_{k=0}^{\infty}\frac{z^k}{k!\Gamma(\rho k+\mu)},
\quad \rho>-1, \quad \mu\in {\mathbb C}.
\end{equation*}
It is easy to see that
\begin{equation} \label{vwf0}
\phi(\rho, \mu; z)\big|_{z=0}=\frac{1}{\Gamma(\mu)}.
\end{equation}
The following  differentiation formulas hold \cite{Wright-1934},
\cite{PsMon}
\begin{equation} \label{dfw}
\frac{d}{dz}\phi(\rho, \mu; z)=\phi(\rho, \mu+\rho; z),
\end{equation}
\begin{equation} \label{dyfw}
D_{0y}^{\nu}y^{\mu-1}\phi(\rho, \mu; -\lambda y^{\rho})=
y^{\rho-\nu-1}\phi(\rho, \rho-\nu;  -\lambda y^{\rho}),
\end{equation}
где $\lambda>0,$ $\rho>-1,$ $\mu, \nu \in {\mathbb R}.$
From the relations (\ref{dfw}) and (\ref{dyfw}) follow
\begin{equation} \label{foeq}
\left(\frac{d}{dz}+\lambda D_{0y}^{\beta}\right)y^{\mu-1}\phi(-\beta, \mu; -\lambda z y^{-\beta})=0,  \quad \beta<1,  \quad \lambda >0.
\end{equation}

For the Wright function following estimates are hold
\cite{PsMon}
\begin{equation} \label{efw0}
|y^{\mu-1}\phi(-\beta, \mu; -\tau y^{-\beta})|
\leq C \tau^{-\theta} y^{\mu+\beta\theta-1}, \quad \tau>0, \quad y>0,
\end{equation}
where $\beta\in(0,1)$ and $\theta\geq 0$ if $\mu\not=0,-1,-2,...,$ and
$\theta\geq -1$ if $\mu=0,-1,-2,...,$ 
\begin{equation} \label{efwi}
|\phi(-\beta, \mu; -z)|\leq C\exp\left(-\sigma z^{\frac{1}{1-\beta}}\right),
\quad z\geq 0,
\end{equation}
where  $\beta\in(0,1),$ $\varepsilon \in {\mathbb R},$
$\sigma<(1-\beta)\beta^{\frac{\beta}{1-\beta}},$ here and below 
$C$ is a pozitive constant.

In paper \cite{Luchko-1999} 
following relation
\begin{equation} \label{iifwn}
\int\limits_0^{\infty}t^n\phi(-\beta, \mu; -t)dt=\frac{n!}{\Gamma(\mu+(n+1)\beta)}, 
\quad n=0,1,...
\end{equation}
was obtained.
In particular for $n=0,$ we have
\begin{equation} \label{iifw}
\int\limits_0^{\infty}\phi(-\beta, \mu; -z)dz=\frac{1}{\Gamma(\mu+\beta)}.
\end{equation}

\section{Wright's function of the matrix argument}

\setcounter{section}{3}
\setcounter{equation}{0}

1. Let $A$ be a square matrix of order $n.$ 
In view of the function $\phi(\rho,\mu;z)$ is analytic everywhere in 
${\mathbb C},$
following series  
\[\phi(\rho,\mu;A)=\sum\limits_{k=0}^{\infty}\frac{A^k }{k!\Gamma(\rho k+\mu)},
\quad \rho>-1, \quad \mu\in {\mathbb C}\]
is convegence for eny  
matrix $A$ given given over the field of complex numbers ${\mathbb C},$ 
and determine the Wright's function of the matrix argument.

Let the matrix $H$ leads the matrix $A$  to Jordan normal form $J(\lambda),$
i.e. 
\begin{equation*}
A=HJ(\lambda)H^{-1},
\end{equation*}
where 
$J(\lambda)={\rm diag}[J_1(\lambda_1),...,J_p(\lambda_p)]$ is the quasydiagonal matrix
with the cells of the form 
\begin{equation*}
J_{k}\equiv J_{k}(\lambda_k)= \left \|
\begin{array}{cccl}
\lambda_k 	& 1   		& \ldots & 0  \\
    		& \lambda_k 	& \ldots & 0 \\
   		&  0  		& \ddots & \vdots \\
   		&     		&        & \lambda_k 
\end{array}
\right \|,  \quad  k=1,...,p,
\end{equation*}
$\lambda_1,...,\lambda_p$ are eigenvalue numbers of the matrix $A,$
and $J_k(\lambda_k)$ are the square matrices of order $r_k+1,$
$\sum\limits_{k=1}^{p}r_k+p=n.$
Then the function $\phi(\rho,\mu;Az)$ can be represented in the form
\begin{equation} \label{wfma}
\phi(\rho,\mu;A z)=H\phi(\rho,\mu;J(\lambda) z)H^{-1},
\end{equation}
where 
\begin{equation*}
\phi(\rho,\mu;J(\lambda) z)={\rm diag}[\phi(\rho,\mu;J_1(\lambda_1) z),...,\phi(\rho,\mu;J_p(\lambda_p) z)],
\end{equation*}
\begin{equation*}
\phi(\rho,\mu; J_{k}(\lambda_k)z)= \left \|
\begin{array}{cccl}
\phi_{\rho,\mu}^{0}(\lambda_k z) 	& \phi_{\rho,\mu}^{1}(\lambda_k z)   	& \ldots 	& \phi_{\rho,\mu}^{r_k}(\lambda_k z)  \\
    					& \phi_{\rho,\mu}^{0}(\lambda_k z) 	& \ldots 	& \phi_{\rho,\mu}^{r_k-1}(\lambda_i z) \\
   					&  0  					& \ddots 	& \vdots \\
   					&     					&        	& \phi_{\rho,\mu}^{0}(\lambda_k z) 
\end{array}
\right \|,
\end{equation*}
\begin{equation*}
\phi_{\rho,\mu}^{m}(\lambda z)=\frac{1}{m!}\frac{\partial^m}{\partial \lambda^m}\phi(\rho,\mu;\lambda z)=\frac{z^m}{m!}\phi(\rho,\mu+\rho m;\lambda z).
\end{equation*}

2.
By using the representation (\ref{wfma}) and the equality (\ref{vwf0}), we obtain
\begin{equation} \label{vwfm0}
\phi(\rho,\mu;A z)\big|_{z=0}=\frac{1}{\Gamma(\mu)}I,
\end{equation}
where $I$ is a identity matrix of order $n.$

3.
Following differentiation  formula holds
\begin{equation} \label{dfwa}
\frac{d}{dz}\phi(\rho, \mu; Az)=A\phi(\rho, \rho+\mu; Az).
\end{equation}
Indeed, by virtue of equality (\ref{dfw}) we get 
\begin{equation*}
\frac{d}{dz}\phi_{\rho,\mu}^{m}(\lambda z)=\lambda\frac{z^m}{m!}\phi(\rho,\mu+\rho +\rho m;\lambda z)+\frac{z^{m-1}}{(m-1)!}\phi(\rho,\mu+\rho m;\lambda z)=
\end{equation*}
\begin{equation*}
=\lambda\phi_{\rho,\mu+\rho}^{m}(\lambda z)+\phi_{\rho,\mu+\rho}^{m-1}(\lambda z).
\end{equation*}
Whence, in turn, we have
\begin{equation} \label{dfwj}
\frac{d}{dz}\phi(\rho, \mu; J(\lambda) z)=J(\lambda)\phi(\rho, \rho+\mu; J(\lambda)z).
\end{equation}
From (\ref{dfwj}), by taking into account following equality  
\begin{equation*}
H\frac{d}{dz}\phi(\rho, \mu; J(\lambda) z)H^{-1}=HJ(\lambda)H^{-1}H\phi(\rho, \rho+\mu; J(\lambda)z)H^{-1}, 
\end{equation*}
we obtain (\ref{dfwa}).

4.
Let all of the eigenvalues of the  matrix $A$ are positive.
Consider the functon
\begin{equation*}
y^{\mu-1}\phi(-\beta,\mu;-A\tau y^{-\beta})=Hy^{\mu-1}\phi(-\beta,\mu;-J(\lambda)\tau y^{-\beta}) H^{-1}.
\end{equation*}
Denoting
\begin{equation*}
w_m\equiv w_m^{\mu}(\tau,y)=\frac{y^{\mu-1}}{m!}\left(-\frac{\tau}{y^{\beta}}\right)^m\phi(-\beta,\mu-m\beta;-\lambda_k\tau y^{-\beta}),
\quad m=0,...,r_k,
\end{equation*}
we can wright 
\begin{equation*}
y^{\mu-1}\phi(-\beta,\mu;-J_k(\lambda_k)\tau y^{-\beta})=
\left \|
\begin{array}{llll}
w_0 	&  w_1  & \ldots &w_{r_k}  \\
	&  w_0  & \ldots &w_{r_k-1} \\
       &     0       & \ddots & \vdots \\
       &             &        &  w_0 
\end{array}
\right \|,
\quad k=1,...,p.
\end{equation*}
In view of (\ref{dyfw}) we get the relation
\begin{equation*}
D_{0y}^{\delta}w_m^{\mu}(\tau,y)=\frac{(-\tau)^m}{m!}y^{\mu-\delta-\beta m-1}\phi(-\beta,\mu-\delta-\beta m;-\lambda\tau y^{-\beta})=w_m^{\mu-\delta}(\tau,y).
\end{equation*}
From this relation, similarly as we obtained the equality (\ref{dfwa}), we obtain the equality
\begin{equation} \label{dyfwa}
D_{0y}^{\delta}y^{\mu-1}\phi(-\beta, \mu; -A\tau y^{-\beta})=y^{\mu-\delta-1}\phi(-\beta, \mu-\delta; -A\tau y^{-\beta}).
\end{equation}

5.
It follows from (\ref{dfwa}) and (\ref{dyfwa})  that
\begin{equation} \label{aphi}
\left(\frac{\partial}{\partial \tau}+AD_{0y}^{\beta}\right)y^{\mu-1}\phi(-\beta,\mu;-A\tau y^{-\beta})=0.
\end{equation}

6.
Let all the eigenvalues of the matrix $A$ be positive,
then the following equality holds
\begin{equation} \label{iifwm}
\int\limits_0^{\infty}\phi(-\beta, \mu; -Az)dz=\frac{1}{\Gamma(\mu+\beta)}A^{-1}.
\end{equation}

Indeed, in view of (\ref{iifwn}) we get
\begin{equation} \label{iifwij}
\int\limits_0^{\infty}\frac{(-z)^m}{m!}\phi(-\beta, \mu-\beta m; -\lambda z)dz=\frac{1}{\Gamma(\mu+\beta)}\frac{(-1)^m}{\lambda^{m+1}}.
\end{equation}
From (\ref{iifwij}) we obtain
$$
\int\limits_0^{\infty}\phi(-\beta, \mu; -J_k z)dz=\frac{1}{\Gamma(\mu+\beta)}
\left\|
\begin{array}{rrrcc}
\frac{1}{\lambda_k} 	& -\frac{1}{\lambda_k^2}   	&\frac{1}{\lambda_k^3}		& \ldots 	& \frac{(-1)^{r_k}}{\lambda_k^{r_k+1}}  \\
    			& \frac{1}{\lambda_k}	 	&-\frac{1}{\lambda_k^2}		& \ldots 	& \frac{(-1)^{r_k-1}}{\lambda_k^{r_k}} \\
   			&  0  				&				& \ddots 	& \vdots \\
   			&     				&        			&		& \frac{1}{\lambda_k} 
\end{array}
\right\|=$$
\begin{equation} \label{iifwj}
=\frac{1}{\Gamma(\mu+\beta)} J_k^{-1}.
\end{equation}
From (\ref{iifwj}), (\ref{wfma}) and equality 
\begin{equation*}
J^{-1}(\lambda)={\rm diag}\left[J_1^{-1}(\lambda_1),...,J_p^{-1}(\lambda_p)\right]
\end{equation*}
follows (\ref{iifwm}).

7.
We denote by $|A(x)|_*$ the scalar function that takes at each point $x$ the largest of the values of the moduli of the elements of the matrix 
$A(x)=\|a_{ij}(x)\|, $ that is $|A(x)|_*=\max\limits_{i,j}|a_{ij}(x)|.$
Similarly, for the vector $b(x)$ with components $b_i(x)$ we denote
$|b(x)|_*=\max\limits_{i}|b_i(x)|.$

From the estimate (\ref{efw0}) follows that
$$|w_m^{\mu}(\tau,y)|\leq \left|\frac{\tau^m}{m!}y^{\mu-\beta m-1}\phi(-\beta,\mu-\beta m;-\lambda\tau y^{-\beta})\right|\leq
$$
\begin{equation*}
\leq C\tau^{m-\theta}y^{\mu-\beta(m-\theta)-1}= C\tau^{-\theta_1}y^{\mu+\beta\theta_1-1},
\end{equation*}
where 
$\theta_1\geq -m,$  при $\mu-\beta m\not= 0,-1,-2,...,$
$\theta_1\geq -m-1,$  при $\mu-\beta m= 0,-1,-2,....$
Thus
\begin{equation} \label{efwa0}
|y^{\mu-1}\phi(-\beta, \mu; -A\tau y^{-\beta})|_*
\leq C \tau^{-\theta} y^{\mu+\beta\theta-1}, \quad \tau>0, \quad y>0,
\end{equation}
где $\beta\in(0,1)$ и $\theta\geq 0$ при $\mu\not=0,-1,-2,...,$ и
$\theta\geq -1$ при $\mu=0,-1,-2,...\,.$ 

8.
From (\ref{efwi}) and (\ref{wfma}) follows the estimate
\begin{equation} \label{efwmai}
|\phi(-\beta, \mu; -Az)|_*\leq C\exp\left(-\sigma z^{\frac{1}{1-\beta}}\right),
\quad z\geq 0,
\end{equation}
where  $\beta\in(0,1),$ $\mu \in {\mathbb R},$
$\sigma<(1-\beta)\left(\lambda\beta^{\beta}\right)^{\frac{1}{1-\beta}},$
$\lambda=\min\limits_{1\leq i\leq p}\{\lambda_i\},$
$\lambda_1, ... , \lambda_p$ are eigenvalues of the matrix $A.$

\section{Main results}

\setcounter{section}{4}
\setcounter{equation}{0}

Consider following function 
$$\Phi_{\alpha}^{\delta}(x)=\int\limits_0^{\infty}e^{B\tau}\prod\limits_{i=1}^m
h_i^{\delta_i}(x_i,\tau)d\tau,
\quad
h_i^{\delta_i}(x_i,\tau)=x_i^{\delta_i-1}\phi(-\alpha_i,\delta_i;-A_i\tau x_i^{-\alpha_i}),$$
where
$\alpha=(\alpha_1,...,\alpha_m),$ $\delta=(\delta_1,...,\delta_m).$

From the estimates (\ref{efwa0}), (\ref{efwmai}) and 
\begin{equation}\label{estexp}
\left|\exp(B\tau)\right|_*\leq C e^{\gamma\tau},
\quad \gamma=\max\limits_{1\leq i\leq q}\{|{\rm Re} \gamma_i|\}.
\end{equation}
where $\gamma_1,...,\gamma_q$ are eigenvalues of matrix $B,$ 
the convergence of the integral 
$\Phi_{\alpha}^{\delta}(x)$ 
follows for all $\alpha_i,\delta_i\in {\mathbb R},$ 
$i=\overline{1,m}.$

A regular solution of system (\ref{bsn}) in domain $\Omega$ is defined as a  
vector function $u=u(x)$ 
satisfying system (\ref{bsn}) at all points $x\in \Omega,$ such that
$D_{0x_i}^{\alpha_i}u \in C(\Omega),$ 
$\prod\limits_{i=1}^mx_i^{1-\mu_i}u(x)\in$ $C(\bar\Omega),$ 
for some positive numbers $\mu_1, \mu_2, ..., \mu_m.$ 


{\bf Theorem 1.}
{\it
Let 
all the eigenvalues of the matrices $A_i$ be positive,
$A_iB=BA_i,$ 
$\mu_{i}<\alpha_i,$ $i=\overline{1,m},$
\begin{equation} \label{busph} 
{\prod\limits_{i=1\atop i\not=j}^m}x_i^{1-\mu_{i}}\varphi_j (x_{(j)})\in 
C(\overline{\Omega^j}),   \quad j=\overline{1,m},
\end{equation}
\begin{equation} \label{busf} 
\prod\limits_{i=1}^mx_i^{1-\mu_i}f(x)\in C(\overline\Omega), 
\end{equation}
$f(x)$ satisfies the Holder condition.
Then there exists a unique regular solution of the problem (\ref{bsn}), (\ref{bu2}), which can be represented as
\begin{equation} \label{bu} 
u(x)=\int\limits_{\Omega_x}G(x-t)f(t)dt+
\sum\limits_{i=1}^m\int\limits_{\Omega^i_x}A_iG(x-t^i)\varphi_i(t_{(i)})dt_{(i)},
\end{equation}
where 
$G(x)=\Phi_{\alpha}^0(x),$
$\Omega_x=\omega_{x_1}\times\omega_{x_2}\times...\times\omega_{x_m},$
$\Omega^i_x=\omega_{x_1}\times...\times\omega_{x_{i-1}}\times\omega_{x_{i+1}}\times...
\times\omega_{x_m},$
$\omega_{x_j}=\{t_j:0<t_j<x_j\},$
$t^i=(t_1,...,t_{i-1},0,t_{i+1},...,t_m).$ 
}

\section{Auxiliary assertions}

\setcounter{section}{5}
\setcounter{equation}{0}

Let us prove some assertions that we need in the proof of Theorem 5.1.

\subsection{Properties of the function  $\Phi_{\alpha}^{\delta}(x)$}

{\bf Lemma 1.} 
{\it The estimate
\begin{equation} \label{boo}
\left|\Phi_{\alpha}^{\delta}(x)\right|_*\leq C \prod\limits_{i=1}^m x_i^{\delta_i+\alpha_i\theta_i-1},
\quad \sum\limits_{i=1}^{m}\theta_i=1,
\end{equation}
holds for all $x_1\in [0,x_{10}],$ 
where 
$\theta_i>0,$ for $\delta_i\not=0,$  and
$\theta_i> -1,$  for $\delta_i=0;$ 
and constant $C$ depends on $x_{10}.$  }

{\bf Proof.}
In view of (\ref{efwa0}) and (\ref{estexp}) we get
$$\left|\Phi_{\alpha}^{\delta}(x)\right|_*\leq C 
\prod\limits_{i=2}^m 
x_i^{\delta_i+\alpha_i\theta_i-1}
\int\limits_0^{\infty}e^{\gamma\tau}\tau^{\theta_1-1}x_1^{\delta_1-1}\phi\left(-\alpha_1,\delta_1;-\tau x_1^{-\alpha_1}\right)d\tau.$$
After replacing $\tau= x_1^{\alpha_1} z,$ we obtain
\begin{equation} \label{boo2}
\left|\Phi_{\alpha}^{\delta}(x)\right|_*\leq C 
\prod\limits_{i=1}^m 
x_i^{\delta_i+\alpha_i\theta_i-1}
\int\limits_0^{\infty}e^{\gamma x_1^{\alpha_1}z}z^{\theta_1-1}\phi(-\alpha_1,\delta_1;-z)dz.
\end{equation}

We represented the integral on the right-hand side (\ref{boo2}) as  following sum
\begin{equation} \label{ej}
J_1(x)+J_2(x)=\left(\int\limits_0^{z_0}+\int\limits_{z_0}^{\infty}\right)z^{-\theta}e^{\gamma x_1^{\alpha_1}z}\phi(-\alpha_1,\mu;-z)dz.
\end{equation}
In view of the boundedness of the function $\phi(-\alpha_1,\mu;-z)$ on any finite interval $[0,z_0],$  we obtain that
\begin{equation} \label{ej1}
|J_1(x)|
\leq
C e^{\gamma x_{10}^{\alpha_1}z_0}\int\limits_0^{z_0}z^{-\theta}dz=C_1e^{\gamma x_{10}^{\alpha_1}z_0}.
\end{equation}
Using the estimate (\ref{efwi}), we have
\[|J_2(x)|\leq C z_0^{-\theta}\int\limits_{z_0}^{\infty}\exp(\gamma x_{10}^{\alpha_1}z-\rho_0 z^{\varepsilon})dz,
\]
where $\rho_0\leq \alpha_1^{\alpha_1/(1-\alpha_1)}(1-\alpha_1),$ $\varepsilon=1/(1-\alpha_1)>1.$
Note that $\gamma x_{10}^{\alpha_1}z-\rho_0 z^{\varepsilon}\leq -z$ for $z\geq z_0=((\gamma x_{10}^{\alpha_1}z+1)/\rho_0)^{(1-\alpha_1)/\alpha_1}>1.$
Therefore
\begin{equation} \label{ej2}
|J_2(x)|\leq C_2 z_0^{-\theta} e^{-z_0}.
\end{equation}
From  (\ref{ej}), (\ref{ej1}), (\ref{ej2}) and (\ref{boo2}) the estimate (\ref{boo}) follows, where
\[C\equiv C(x_{10})=C_1e^{\gamma x_{10}^{\alpha_1}z_0}+C_2 z_0^{-\theta} e^{-z_0}.\]
Lemma 1 is proved.

{\bf Remark.}
From the equalities (\ref{dfwa}) and (\ref{dyfwa}), it is easy to see that formula
\begin{equation} \label{dpad}
D_{0x_i}^{\nu}\Phi_{\alpha}^{\delta}(x)=
\Phi_{\alpha}^{(\delta_1,...,\delta_i-\nu,...,\delta_m}(x)
\end{equation}
is valid for $\nu=0$ and $\nu\in{\mathbb N}.$
In other cases it is required that the function $\Phi_{\alpha}^{\delta}(x)$ has an integrable singularity for $x_i=0.$
As follows from the estimate (\ref{boo}), for this is sufficient $\delta_i+\alpha_i> 0.$

{\bf Lemma 2.} 
{\it Let $A_iB=BA_i,$ 
$\delta_i+\alpha_i>0$  $(i=\overline{1,m}),$
then the equality
\begin{equation}\label{lpabmn}
\left(\sum\limits_{i=1}^{m}A_iD_{0x_i}^{\alpha_i} -B\right)\Phi_{\alpha}^{\delta}(x)=\prod\limits_{i=1}^{m}\frac{x_i^{\delta_i-1}}{\Gamma(\delta_i)}I,
\end{equation}
holds for all $x\in \Omega,$ where $I$ is identity matrix.}

{\bf Proof.}
By virtue of (\ref{dyfwa}) we obtain the expression
\begin{equation}\label{bxkh}
D_{0x_k}^{\alpha_k}\Phi_{\alpha}^{\delta}(x)=
\int\limits_0^{\infty}e^{B\tau} D_{0x_k}^{\alpha_k} h_k^{\delta_k}(x_k,\tau)\prod\limits_{i=1 \atop i\not=k}^{m}h_i^{\delta_i}(x_i,\tau)d\tau, 
\quad k=1,...,m
\end{equation}
 for the derivatives with respect to $x_k.$ We transform the previous formula for $k=1.$
Using the equation (\ref{aphi}), by the formula of integration by parts, and the relations (\ref{vwfm0}) and (\ref{efwmai}), we obtain
$$D_{0x_1}^{\alpha_1}\Phi_{\alpha}^{\delta}(x)=
-\int\limits_0^{\infty}e^{B\tau} \left[\frac{\partial}{\partial\tau}h_1^{\delta_1}(x_1,\tau)\right]\prod\limits_{i=2}^{m}h_i^{\delta_i}(x_i,\tau)d\tau= $$
$$=-e^{B\tau}\prod\limits_{i=1}^{m}h_i^{\delta_i}(x_i,\tau)\Big|_{\tau=0}^{\tau=\infty}+
\int\limits_0^{\infty} h_1^{\delta_1}(x_1,\tau)\frac{\partial}{\partial\tau}\left[e^{B\tau}\prod\limits_{i=2}^{m}h_i^{\delta_i}(x_i,\tau)\right]d\tau=$$
$$=\prod\limits_{i=1}^{m}\frac{x_i^{\delta_i-1}}{\Gamma(\delta_i)}I+B\Phi_{\alpha}^{\delta}(x)+
\int\limits_0^{\infty} h_1^{\delta_1}(x_1,\tau)e^{B\tau}\frac{\partial}{\partial\tau}\prod\limits_{i=2}^{m}h_i^{\delta_i}(x_i,\tau)d\tau=$$
\begin{equation}\label{bxkh1}
=\prod\limits_{i=1}^{m}\frac{x_i^{\delta_i-1}}{\Gamma(\delta_i)}I+B\Phi_{\alpha}^{\delta}(x)+
\sum\limits_{k=2}^{m}\int\limits_0^{\infty} e^{B\tau}\left[\frac{\partial}{\partial\tau}h_k^{\delta_k}(x_k,\tau)\right]\prod\limits_{i=1 \atop i\not=k}^{m}h_i^{\delta_i}(x_i,\tau)d\tau.
\end{equation}
From the relations (\ref{bxkh}) and (\ref{bxkh1}) follows that
\begin{equation*}
\left(\sum\limits_{k=1}^{m}A_kD_{0x_k}^{\alpha_k} -B\right)\Phi_{\alpha}^{\delta}(x)=\prod\limits_{i=1}^{m}\frac{x_i^{\delta_i-1}}{\Gamma(\delta_i)}I+
\end{equation*}
\begin{equation*}
+\sum\limits_{k=2}^{m}\int\limits_0^{\infty} e^{B\tau}\prod\limits_{i=1 \atop i\not=k}^{m}h_i^{\delta_i}(x_i,\tau)\left[\frac{\partial}{\partial\tau}+
A_kD_{0x_k}^{\alpha_k}\right]h_k^{\delta_k}(x_k,\tau)d\tau.
\end{equation*}
From this, taking into account (\ref{aphi}), we obtain (\ref{lpabmn}).
Lemma 2 is proved.

\subsection{Representation of the solution}

{\bf Lemma 3.} 
{\it Let $A_iB=BA_i,$ then any regular in the domain $\Omega$ solution $u(x)$ 
of the problem (\ref{bsn}) - (\ref{bu2})  can be represented as (\ref{bu}).}

{\bf Proof.}
Let the matrix $V(x)$ be a solution of equation
\begin{equation} \label{bsnv}
\sum\limits_{i=1}^{m}D_{0x_i}^{\alpha_i} V(x)A_i=V(x)B+I,
\end{equation}
satisfying conditions
\begin{equation} \label{busnv}
V(x^i)=0,\quad 
i=\overline{1,n},
\end{equation}
where $I$ is the identity matrix of order $n.$
In view of Lemma 2 the function
$$V(x)=\int\limits_{\Omega_x}G(t)dt=\Phi_{\alpha}^{(1,...,1)}(x),$$
is a solution of equation
\begin{equation} \label{bsnv1}
\sum\limits_{i=1}^{m}A_iD_{0x_i}^{\alpha_i} V(x)=BV(x)+I.
\end{equation}
Hence, taking into account the fact that the matrix 
$\Phi_{\alpha}^{\delta}(x)$ commutes with the matrices $A_i$ and $B,$ 
we obtain that $V(x)$ is a solution of the equation (\ref{bsnv}).

It follows from (\ref{boo}) that the estimate $$|V(x)|_* \leq C \prod\limits_{i=1}^{m}x_i^{\alpha_i\theta_i}, \quad \theta_i> 0,$$ 
from which, in turn, follows (\ref{busnv}).
That is, $V(x)$ is the solution of the problem (\ref{bsnv}), (\ref{busnv}).

By virtue of the integration by parts formula and (\ref{busnv}), we obtain
$$\int\limits_{\Omega_x\setminus\Omega_{\varepsilon}}V(x-t)
\sum\limits_{i=1}^{m}A_iD_{0t_i}^{\alpha_i}u(t)dt=-\sum\limits_{i=1}^{m}
\int\limits_{\Omega_x\setminus\Omega_{\varepsilon}}
\left[\frac{\partial}{\partial t_i}V(x-t)\right]A_iD_{0t_i}^{\alpha_i-1}u(t)dt-$$
$$-\sum\limits_{i=1}^{m}\int\limits_{\Omega_x^i\setminus\Omega^i_{\varepsilon}}
V(x-t)\big|_{t_i=\varepsilon_i}A_iD_{0\varepsilon_i}^{\alpha_i-1}u(t)dt_{(i)},
\quad i=1,...,m,$$
where 
$\Omega_{\varepsilon}=\omega_{\varepsilon_1}\times...\times\omega_{\varepsilon_n},$
$\Omega^i_{\varepsilon}=\omega_{\varepsilon_1}\times...
\times\omega_{\varepsilon_{i-1}}\times\omega_{\varepsilon_{i+1}}\times...
\times\omega_{\varepsilon_n},$
$\omega_{\varepsilon_j}=\{t_j:0<t_j<\varepsilon_j\},$
$\varepsilon=(\varepsilon_1,...,\varepsilon_n).$
Passing in the last equality to the limit as $\varepsilon_i \to 0,$
taking into account (\ref{bsn}), (\ref{bu2}), (\ref{bsnv}) and equality \cite[p. 34]{Nakhushev-2003}
$$\int\limits_0^xv_1(t)D_{0t}^{\nu}v_2(t)dt=\int\limits_0^x\left[D_{xt}^{\nu}v_1(t)\right]v_2(t)dt,
\quad \nu<0,$$
we obtain
\begin{equation} \label{biuxy}
\int\limits_{\Omega_x}u(t)dt=
\int\limits_{\Omega_x}V(x-t)f(t)dt+
\sum\limits_{i=1}^{m}\int\limits_{\Omega_x^i}V(x-t^i)A_i\varphi_i(t_{(i)})dt_{(i)}.
\end{equation}
Differentiating the equality (\ref{biuxy}) over all $x_i$ and taken into account (\ref{busnv}), we obtain
\begin{equation} \label{buxy}
u(x)=
\int\limits_{\Omega_x}V_{x_1...x_m}(x-t)f(t)dt+
\sum\limits_{i=1}^{m}\int\limits_{\Omega_x^i}V_{x_1...x_m}(x-t^i)A_i\varphi_i(t_{(i)})dt_{(i)}.
\end{equation}
From (\ref{buxy}), and the equality $V_{x_1 ... x_m}(x)=G(x),$ we get (\ref{bu}).
Lemma 3 is proved.

\subsection{Properties of the fundamental solution}

{\bf Lemma 4.}   
{\it The following estimates hold
\begin{equation} \label{bogab}
|G(x)|_*\leq C \prod\limits_{i=1}^m x_i^{\alpha_i\theta_i-1}, 
\quad \sum_{i=1}^{m}\theta_i=1, \quad  \theta_i> -1;
\end{equation}
\begin{equation} \label{bga00}
|D_{0x_s}^{\alpha_s-1}G(x)|_*\leq C x_s^{1-\alpha_s}\prod\limits_{i=1}^m x_i^{\alpha_i\theta_i-1},
\quad \sum_{i=1}^{m}\theta_i=1, \quad  
\theta_i> \left\{
\begin{array}{l}
0, i=s,  \\ 
-1, i\not=s;
\end{array}
\right. 
\end{equation}
\begin{equation}\label{bf1511}
|D_{0x_s}^{\alpha_s}G(x)|_*\leq C x_s^{-\alpha_s}
\prod\limits_{i=1}^m x_i^{\alpha_i\theta_i-1},
\quad \sum_{i=1}^{m}\theta_i=1, \quad  \theta_i> -1,
\end{equation}
where $C$ is the positive constant.}

{\bf Lemma 5.} 
{\it The equality
$$\sum\limits_{i=1}^{m}D_{0x_i}^{\alpha_i} G(x)=BG(x).$$
holds for all $x\in \Omega.$}

Lemmas 4 and 5 follow from Lemmas 1 and 2 and the formula (\ref{dpad}).

{\bf  Lemma 6.} 
{\it Let the conditions of Theorem 1 be fulfilled, then the function $u(x)$ defined by the equality 
(\ref{bu}) is a solution of the equation (\ref{bsn}),
such that} $D_{0x_i}^{\alpha_i}u\in C(\Omega).$
                
{\bf Proof.}
Denote 
$$u_f(x)=\int\limits_{\Omega_x}G(x-t)f(t)dt, \quad 
u^i(x)=A_i\int\limits_{\Omega^i_x}G(x-t^i)\varphi_i(t_{(i)})dt_{(i)}, \quad i=1,...,m.$$ 
From (\ref{bf1511}) we obtain 
\begin{equation} \label{bttt1}
|D_{x_st_s}^{\alpha_s}G(x-t^i)|_*<C x_i^{\alpha_i\theta_i-1}(x_s-t_s)^{-\alpha_s}
\prod\limits_{j=1 \atop j\not=i}(x_j-t_j)^{\alpha_j\theta_j-1}, \quad s\not=i,
\end{equation}
\begin{equation} \label{bttt2}
|D_{0t_s}^{\alpha_s}G(x-t^s)|_*<C x_s^{\alpha_s\theta_s-\alpha_s-1}
\prod\limits_{j=1 \atop j\not=i}(x_j-t_j)^{\alpha_j\theta_j-1},
\end{equation}
where 
$\theta_i> \left\{
\begin{array}{l}
1, i=s,  \\ 
-1, i\not=s,
\end{array}
\right. 
\quad \sum\limits_{i=1}^{n}\theta_i=1.$

By virtue of formula \cite[c. 99]{Podlubny}
$$D_{0x}^{\nu}\int\limits_0^xv_1(x-t)v_2(t)dt=\int\limits_0^xv_1(x-t)D_{0t}^{\nu}v_2(t)dt+$$
$$+v_1(x)\lim\limits_{x\rightarrow 0}D_{0x}^{\nu-1}v_2(t), \quad 0<\nu<1, 
$$
and estimates (\ref{bttt1}) and (\ref{bttt2}), the inclusions
\begin{equation} \label{bi}
D_{0x_s}^{\alpha_s}u^i(x)=
A_i\int\limits_{\Omega^i_x}[D_{x_st_s}^{\alpha_s}G(x-t^i)]\varphi_i(t_{(i)})dt_{(i)}
\in C(\Omega), \quad i\not=s,
\end{equation}
\begin{equation} \label{bii}
D_{0x_s}^{\alpha_s}u^s(x)=
A_s\int\limits_{\Omega^s_x}[D_{0x_s}^{\alpha_s}G(x-t^s)]\varphi_s(t_{(s)})dt_{(s)}
\in C(\Omega). 
\end{equation}
and equalities
\begin{equation} \label{biii}
\sum\limits_{j=1}^m A_jD_{0x_j}^{\alpha_j} u^i(x)=Bu_i(x), \quad i=1,...,m
\end{equation}
are valid.

Consider following integrals
$$J(x, t_{(k)})=\int\limits_0^{x_k}[D_{x_k t_k}^{\alpha_k-1}G(x-t)]f(t)dt_k,$$
$$J_{\varepsilon}(x, t_{(k)})=\int\limits_0^{x_k-\varepsilon}
[D_{x_k t_k}^{\alpha_k-1}G(x-t)]f(t)dt_k.$$
It is obvious that
$\lim\limits_{\varepsilon\rightarrow 0}J_{\varepsilon}(x, t_{(k)})=J(x, t_{(k)}).$
By (\ref{bf1511}), and the fact that $f(x)$ satisfies the Holder condition, we obtain
$$\left|D_{x_k t_k}^{\alpha_k}G(x-t)[f(t_k,t_{(k)})-f(x_k,t_{(k)})]\right|_*\leq$$
$$\leq C(x_k-t_k)^{q+\alpha_k(\theta_k-1)-1}
\prod\limits_{j=1 \atop j\not=k}^nt_j^{\mu_j-1}(x_j-t_j)^{\alpha_j\theta_j-1},$$
here $\theta_k>1-\frac{q}{\alpha_k},$ $\sum\limits_{j=1}^{n}\theta_j=1.$
Hence it is easy to see that the integral on the right-hand side of 
$$\frac{\partial}{\partial x_k}J_{\varepsilon}(x, t_{(k)})=
\int\limits_0^{x_k-\varepsilon}
\left[\frac{\partial}{\partial x_k}D_{x_kt_k}^{\alpha_k-1}G(x-t)\right]
[f(t_k,t_{(k)})-f(x_k,t_{(k)})]dt_k -$$
$$-D_{x_kt_k}^{\alpha_k-1}G(x-t)|_{t_k=0}^{t_k=x_k-\varepsilon}f(x_k,t_{(k)})+
[D_{x_k,x_k-\varepsilon }^{\alpha_k-1}G(x-t)]f(x_k-\varepsilon,t_{(k)})$$
converges uniformly on the set
$\Omega\times\Omega_x^k$ for all $q\in(0,1].$
Therefore
$$\lim\limits_{\varepsilon \rightarrow 0}
\frac{\partial}{\partial x_k}J_{\varepsilon}(x, t_{(k)})=
\frac{\partial}{\partial x_k}J(x, t_{(k)})=
\left[D_{x_kt_k}^{\alpha_k-1}G(x-t)|_{t_k=0}\right]f(x_k,t_{(k)}) +$$
$$+\int\limits_0^{x_k}[D_{x_kt_k}^{\alpha_k}G(x-t)][f(t_k,t_{(k)})-f(x_k,t_{(k)})]dt_k.$$
From the latter, taking into account
$$\Big|\frac{\partial}{\partial x_k}J_{\varepsilon}(x, t_{(k)})\Big|_*
\leq C x_k^{q-\alpha_k+\alpha_k\theta_k}
\prod\limits_{j=1 \atop j\not=k}^nt_j^{\mu_j-1}(x_j-t_j)^{\alpha_j\theta_j-1},
$$
we get
\begin{equation} \label{biv}
D_{0x_k}^{\alpha_k}u_f(x)=\frac{\partial}{\partial x_k}\int\limits_{\Omega_x}
[D_{x_kt_k}^{\alpha_k-1}G(x-t)]f(t)dt=
\int\limits_{\Omega_x}\frac{\partial}{\partial x_k}J(x, t_{(k)})dt_{(k)}\in C(\Omega).
\end{equation}
In veiw of (\ref{bsnv1}), we have
$$\int\limits_{\Omega_x}\left[\sum\limits_{k=1}^{m}A_kD_{0t_k}^{\alpha_k}-B\right]u_f(t)dt=
\int\limits_{\Omega_x}\left[\sum\limits_{k=1}^{m}A_kD_{x_kt_k}^{\alpha_k}-B\right]V(x-t)f(t)dt=
\int\limits_{\Omega_x}f(t)dt.$$
Which implies that
\begin{equation} \label{bv}
\sum\limits_{k=1}^{m}A_kD_{0x_k}^{\alpha_k}u_f(x)-Bu_f(x)=f(x).
\end{equation}
The validity of the lemma 6 follows from (\ref{bi}) - (\ref{bv}).
Lemma 6 is proved.

{\bf Lemma 7.} 
{\it 
Let the function $\varphi_j(x_{(j)})$  satisfies the condition (\ref{busph}), 
then the following relations
\begin{equation} \label{blxp}
\lim\limits_{x_s\rightarrow 0}D_{0x_s}^{\alpha_s-1}\int\limits_{\Omega_x^j} 
G(x-t^j)\varphi_j(t_{(j)})dt_{(j)}=0,
\quad s\not=j, \quad
x_{(s)}\in\Omega^s\setminus\Omega_{\varepsilon}^s,
\end{equation}
\begin{equation} \label{blxg}
\lim\limits_{x_s\rightarrow 0}D_{0x_s}^{\alpha_s-1}\int\limits_{\Omega_x^s} 
G(x-t^s)\varphi_s(t_{(s)})dt_{(s)}=\varphi_s(x_{(s)}),
\quad
x_{(s)}\in\Omega^s\setminus\Omega_{\varepsilon}^s
\end{equation}
hold, 
and the limits are uniform on any closed subset of the domain $\Omega^s.$ }

{\bf Proof.}
In view of (\ref{bga00}) we have the estimate
\begin{equation} \label{blxo1}
\left|D_{0x_s}^{\alpha_s-1}u^i(x)\right|_*\leq C x_s^{\alpha_s\theta_s+\mu_{sj}-\alpha_s}, 
\quad s\not=j, 
\end{equation}
for $x_{(s)}\in\Omega^s\setminus\Omega_{\varepsilon}^s.$ 
From this estimate follows (\ref{blxp}).

Consider the integral
$$D_{0x_s}^{\alpha_s-1}u^s(x)=
A_sD_{0x_s}^{\alpha_s-1}\int\limits_{\Omega_x^s} G(x-t^s)\varphi_s(t_{(s)})dt_{(s)}=$$
\begin{equation} \label{bdus}
=A_s\left(\int\limits_{\Omega_{\varepsilon}^s} +\int\limits_{\Omega_x^s\setminus\Omega_{\varepsilon}^s} \right)
D_{0x_s}^{\alpha_s-1}G(x_s,t_{(s)})
\varphi_s(x_{(s)}-t_{(s)})dt_{(s)},
\end{equation}
where
$G(x_s,t_{(s)})=G(t)|_{t_s=x_s},$ 
$\Omega^i_{\varepsilon}=\omega_{\varepsilon_1}\times...
\times\omega_{\varepsilon_{i-1}}\times\omega_{\varepsilon_{i+1}}\times...
\times\omega_{\varepsilon_n},$
$\omega_{\varepsilon_j}=\{t_j:0<t_j<\varepsilon_j\},$
$\varepsilon=(\varepsilon_1,...,\varepsilon_n).$
The limit of the second integral for $x_s\rightarrow 0$ is zero, due to the estimate
$$\left|D_{0x_s}^{\alpha_s-1}G(x_s,t_{(s)})\right|\leq C x_s^{\alpha_s\theta}, 
\quad  0<\theta\leq 1,
\quad  x_{(s)}\in\Omega^s\setminus\Omega_{\varepsilon}^s,
$$
and the boundedness of the integral
$\int\limits_{\Omega_x^s\setminus\Omega_{\varepsilon}^s}\varphi_s(x_{(s)}-t_{(s)})dt_{(s)}.$  
We denote the first integral  $I_1(x),$ then
$$I_1(x)=
A_s\int\limits_{\Omega_{\varepsilon}^s}D_{0x_s}^{\alpha_s-1}G(x_s,t_{(s)}) 
\left[\varphi_s(x_{(s)}-t_{(s)})-\varphi_s(x_{(s)})\right]dt_{(s)}+
$$
$$+A_s\left[\int\limits_{\Omega_{\varepsilon}^s}D_{0x_s}^{\alpha_s-1}G(x_s,t_{(s)})dt_{(s)}\right] 
\varphi_s(x_{(s)})=
$$
$$=
A_s\int\limits_{\Omega_{\varepsilon}^s}
\left[\int\limits_0^{\infty}e^{B \tau}h_s^{1-\alpha_s}(x_s,\tau)
\prod\limits_{i=1 \atop i\not=s}^mh_i^0(t_i,\tau)d\tau\right]
\left[\varphi_s(x_{(s)}-t_{(s)})-\varphi_s(x_{(s)})\right]dt_{(s)}+
$$
\begin{equation} \label{bdusi1}
+A_s\left[
\int\limits_0^{\infty}e^{B \tau}h_s^{1-\alpha_s}(x_s,\tau)d\tau
\int\limits_{\Omega_{\varepsilon}^s}\prod\limits_{i=1 \atop i\not=s}^mh_i^0(t_i,\tau)dt_{(s)}\right] 
\varphi_s(x_{(s)}).
\end{equation}
Using the fact that by virtue of (\ref{dyfwa}) 
$$\int\limits_0^{\varepsilon_i}\frac{1}{t_i}
\phi\left(-\alpha_i,0;-A_i\tau t_i^{-\alpha_i}\right)dt_i=
\phi\left(-\alpha_i,1;-A_i\tau \varepsilon_i^{-\alpha_i}\right),$$ 
and then replacing the integration variable, we transform the integral
$$
\int\limits_0^{\infty}e^{B \tau}h_s^{1-\alpha_s}(x_s,\tau)d\tau
\int\limits_{\Omega_{\varepsilon}^s}\prod\limits_{i=1 \atop i\not=s}^mh_i^0(t_i,\tau)dt_{(s)}=
\int\limits_0^{\infty}e^{B \tau}h_s^{1-\alpha_s}(x_s,\tau)
\prod\limits_{i=1 \atop i\not=s}^mh_i^1(\varepsilon_i,\tau)d\tau=$$
\begin{equation} \label{bfyz}
=\int\limits_0^{\infty}
\phi(-\alpha_s,1-\alpha_s;-A_sz) F(x_s,z)dz.
\end{equation}
where
$F(x_s,z)=e^{B x_s^{\alpha_s}z}
\prod\limits_{i=1 \atop i\not=s}^n
\phi\left(-\alpha_i,1;-A_i\varepsilon_i^{-\alpha_i}x_s^{\alpha_s}z\right).$

It follows from the estimate (\ref{efw0}) that there exists a uniform limit for all $z\in[0,z_0]$ $z_0<\infty$ 
\begin{equation} \label{blxsf}
\lim\limits_{x_s\rightarrow 0}F(x_s,z)=I,
\end{equation}
and, that $|F(x_s,z)|_*\leq \exp(\gamma x_{s0}^{\alpha_s} z),$ for each finite $x_s\leq x_{s0}.$ 
It follows from the latter that the integral (\ref{bfyz}) converges uniformly in all $x_s\in[0, x_{s0}].$
Passing to the limit in the integral (\ref{bfyz}) as $x_s\rightarrow 0,$ by
taking into account (\ref{blxsf}) and the formula (\ref{iifwm}), we get
$$\lim\limits_{x_s\rightarrow 0}
\int\limits_0^{\infty}e^{B \tau}h_s^{1-\alpha_s}(x_s,\tau)d\tau
\int\limits_{\Omega_{\varepsilon}^s}\prod\limits_{i=1 \atop i\not=s}^mh_i^0(t_i,\tau)dt_{(s)}=
$$
\begin{equation} \label{bi11}
=\int\limits_0^{\infty}\phi(-\alpha_s,1-\alpha_s;-A_sz)dz=A_s^{-1}. 
\end{equation}
The function $\varphi_s(t_{(s)})$ is continuous on $[x-\varepsilon,x],$  therefore 
$$\omega(\varepsilon)=\sup |\varphi_s(x_{(s)}-t_{(s)})-\varphi_s(x_{(s)})|\rightarrow 0$$ 
for $\varepsilon\rightarrow 0.$ 
Because of the arbitrariness of the choice of $\varepsilon$
and (\ref{bi11}),
for $x_s\rightarrow 0$ the first term in (\ref{bdusi1}) tends to zero, and the second to $\varphi_s(x _{(s)}).$
Thus
$$\lim\limits_{x_s\rightarrow 0}I_1(x)=\varphi_s(x_{(s)}).$$
From the latter together with (\ref{bdus}) it follows (\ref{blxg}).
Lemma 7 is proved.

\section{Proof of Theorem 1}
\setcounter{section}{6}
\setcounter{equation}{0}

{\bf Proof.}
Taking into account that  $f^*(x)=\prod\limits_{i=1}^mx_i^{1-\mu_i}f(x)\in C(\overline{\Omega})$ 
and using  (\ref{busf})  and (\ref{bogab}) we get 
\begin{equation} \label{bufo}
|u_f(x)|_*\leq C\prod\limits_{i=1}^m x_i^{\mu_i+\alpha_i\theta_i-1}, 
\quad \theta_i>0, \quad \sum\limits_{i=1}^m \theta_i=1.
\end{equation} 
It follow from  
(\ref{bufo}) that 
$$\lim\limits_{x_s\rightarrow 0}D_{0x_s}^{\alpha_s-1}u_f(x)=0.$$
From Lemma 7 and the last relation, the fulfillment of the boundary conditions (\ref{bu2}) follows.
From the estimate (\ref{bufo}) it also follows that
$\prod\limits_{i=1}^n x_i^{1-\mu_i}u_f\in C(\overline{\Omega}).$

Using (\ref{busph}) and (\ref{bogab}) give 
$$|u_k(x)|_*\leq 
Cx_k^{\alpha_k\theta_{k}-1}
\prod\limits_{i=1 \atop i\not= k}^m x_i^{\alpha_i\theta_{i}+\mu_{i}-1}.$$ 
From the last inequality it follows that
$\prod\limits_{i=1}^m x_i^{1-\mu_i}(u-u_f)\in C(\overline{\Omega}).$
The foregoing, together with Lemma 6
proves the existence of a regular solution of the problem
(\ref{bsn}) -- (\ref{bu2}).
The uniqueness of the solution of the problem follows from Lemma 5. 
Theorem 1 is proved.






 \bigskip \smallskip

 \it

\noindent
Institute of Applied Mathematics and Automation of KBSC RAS\\[4pt]
"Shortanov" $\,$ Str., 89 A, Nal'chik -- 360000, RUSSIA  \\[4pt]
e-mail: mamchuev@rambler.ru


\begin{thebibliography}{99}
\normalsize




\bibitem{Nakhushev-2003} 
Nakhushev A.~M. {\it Fractional calculus and its applications}, Moscow: Fizmatlit, 2003. 
(In Russian). 

\bibitem{Clement-2000} 
Clement Ph., Gripenberg G., Londen S-O. {\it Schauder estimates for equations with fractional derivatives},   
Trans. of the Amer. Math. Soc., {\bf 352}:5, (2000), 2239--2260.  

\bibitem{Wright-1933}  
Wright~E.~M. 
{\it On the coefficients of power series having exponential singularities}, 
J. London Math. Soc.,  {\bf 8:}29, (1933),  71--79.  

\bibitem{Wright-1934}  
Wright~E.~M. 
{\it The asymptotic expansion of the generalized Bessel function}, 
Proc. London Math. Soc. Ser. II,  {\bf 38}, 257--270. 

\bibitem{Clement-1999} 
Clement Ph., Gripenberg G., Londen S-O. {\it Holder regularity for a linear fractional evolution equation}, 
Progr. Nonlinear Differ. Equat. and Their Appl., {\bf 35}, (1999), 62--82. 


\bibitem{Pskhu-2003} 
Pskhu A.~V. {\it Solution of a boundary value problem for a fractional partial differential equation}, 
Differential Equation, {\bf 39:}8,  (2003),  1092-1099.
(In Russian). 

\bibitem{PsMon} 
Pskhu A.~V. {\it Fractional partial differential equations},  Moscow: Nauka, 2005. (In Russian).


\bibitem{Mamchuev-2009-1} 
Mamchuev M.~O.  {\it A boundary value problem for a first-order equation with a partial derivative of a fractional order with variable coefficients},
Reports of Circassian International Academy of Sciences. {\bf 11:}1, (2009), 32--35. (In Russian).

\bibitem{Mamchuev-2009-2}  
Mamchuev M.~O.  {\it Cauchy problem in non-local statement for first order equation with partial derivatives of fractional order with variable coefficients},
Reports of Circassian International Academy of Sciences. {\bf 11:}2, (2009), 21--24. (In Russian).

\bibitem{MamMon} 
Mamchuev M.~O. 
{\it Boundary value problems for equations and systems with the partial derivatives of fractional order}, Nalchik: Publishing house KBSC of RAS, 2013. (In Russian).

\bibitem{Pskhu-2011} 
Pskhu A.~V.
{\it Boundary value problem for a multidimensional fractional partial differential equation}, 
Differential Equation, {\bf 47:}3,  (2011),  385-395.
(In Russian).

\bibitem{Mamchuev-2008} 
Mamchuev M.~O.  
{\it Boundary value problem for a system of fractional partial differential equations}, 
Differential Equations, {\bf 44:}12, (2008), 1737--1749.   


\bibitem{Mamchuev-2017} 
Mamchuev M.~O.
{\it Boundary value problem for a linear system of equations with the partial derivatives of fractional order}, 
Chelyabinsk Physical and Mathematical Journal, {\bf 2:}3, (2017), 295--311. (In Russian).  


\bibitem{Luchko-1999}  
Gorenflo R., Luchko Yu., Mainardi F.  
{\it Analytical properties and applications of the Wright function}, 
Fractional Calculus and Applied Analysis,    {\bf 2:}4, (1999),  383--414.  

\bibitem{Podlubny} 
{Podlubny I. }
{\it Fractional differential equations}, New-York: Acad. press, 1999. 










\end{thebibliography}
\end{document}